\documentclass[12pt,twoside,reqno]{amsart}
\usepackage{amsxtra,amscd}
\usepackage{graphicx}
\usepackage{amsmath}
\usepackage{amsfonts}
\usepackage{amssymb}

\newtheorem{theorem}{Theorem}[section]
\newtheorem{lemma}[theorem]{Lemma}
\theoremstyle{definition}
\newtheorem{definition}[theorem]{Definition}

\newtheorem{proposition}[theorem]{Proposition}
\newtheorem{remark}[theorem]{Remark}

\newtheorem{corollary}[theorem]{Corollary}

\numberwithin{equation}{section}

\begin{document}
\title{The Conjugates of Algebraic Schemes}

\author{Feng-Wen An}
\address{School of Mathematics and Statistics, Wuhan University, Wuhan,
Hubei 430072, People's Republic of China} \email{fwan@amss.ac.cn}


\begin{abstract}
Fixed an algebraic scheme $Y$. We suggest a definition for the
conjugate of an algebraic scheme $X$ over $Y$ in an evident manner;
then $X$ is said to be Galois closed over $Y$ if $X$ has a unique
conjugate over $Y$. Now let $X$ and $Y$ both be integral and let $X$
be Galois closed over $Y$ by a surjective morphism $\phi$ of finite
type. Then $\phi^{\sharp}(k(Y))$ is a subfield of $k(X)$ by $\phi$.
The main theorem of this paper says that $k\left( X\right)
/\phi^{\sharp}(k\left( Y\right)) $ is a Galois extension and the
Galois group $Gal\left(k(X)/\phi^{\sharp}(k\left( Y\right)) \right)
$ is isomorphic to the group of $k-$automorphisms of $X$ over $Y$ if
$k\left( X\right) /\phi^{\sharp}(k(Y))$ is separably generated.
\end{abstract}

\maketitle





\bigskip

\section*{\textbf{\textrm{Contents}}}

\qquad \textsf{Introduction}

\qquad\quad \textsf{0.1 Background and Motivation}

\qquad\quad \textsf{0.2 Main Theorem of the Paper, An Introduction}

\qquad\quad \textsf{0.3 Outline of the Proof for the Main Theorem}

\qquad \textsf{1. Definition for Galois Closed Schemes}

\qquad\quad \textsf{1.1 Galois Closed Schemes}

\qquad\quad \textsf{1.2 Galois Closures}

\qquad \textsf{2. Statement of The Main Theorem}

\qquad\quad \textsf{2.1 Notation}

\qquad\quad \textsf{2.2 Statement of The Main Theorem}

\qquad \textsf{3. Proof of The Main Theorem}

\qquad\quad \textsf{3.1 Definition for Conjugations of a Field}

\qquad\quad \textsf{3.2 Definition for Conjugations of an Open Set}

\qquad\quad \textsf{3.3 Conjugations and Galois Extensions}

\qquad\quad \textsf{3.4 Preliminaries: $k-$Affine Structures}

\qquad\quad \textsf{3.5 Conjugations and Galois Closed Varieties}

\qquad\quad \textsf{3.6 Proof of the Main Theorem}

\qquad \textsf{References}

\bigskip

\bigskip

\section*{Introduction}

\bigskip

\subsection{Background and Motivation}

Let us begin with Weil's notion$^{[9]}$, the conjugate $X^{\sigma}$
of a classical variety $X$ defined over an algebraic extension of a
field $k$, where $\sigma$ is a given $k-$automorphism of
$\overline{k}$. $X$ and $X^{\sigma}$ behaving like conjugates of a
field, are almost of the same algebraic properties; however, their
topological properties are very different from each other in
general.

For example, if $k$ is a number field, Serre shows us an
example$^{\left[ 8\right] }$ that $X^{an}$ and $\left(
X^{\sigma}\right) ^{an}$ can not be topologically equivalent spaces,
where $X^{an}$ denotes the analytical space associated with $X$.
There can be the relation $\pi_{1}^{top}\left( X\right)
\ncong\pi_{1}^{top}\left( X^{\sigma}\right)$ for a variety $X$ and a
$ \sigma\in Aut_{k}\left( \overline{k}\right) $. By a
theorem$^{\left[ 3\right] }$ of Grothendieck, the profinite
completions of their topological fundamental groups are equal, that
is, $\pi_{1}^{alg}\left( X\right) \cong\pi_{1}^{alg}\left(
X^{\sigma}\right) .$ It has been still unknown why there exists such
a phenomenon at least to the author's knowledge.

On the other hand, there have been various discussions$^{[4-7]}$
which use the relevant data of varieties to describe class field
theory, especially use fundamental groups to describe nonabelian
theory in recent years.

In this paper we will extend in an obvious manner the conjugate of a
classical variety to a more general case, the conjugate of an
algebraic scheme. Then an algebraic scheme is said to be Galois
closed if it has a unique conjugate. We try to use these relevant
data of such schemes to obtain some information of a Galois
extension of a field. With Galois closed schemes, we believe that in
the future we can obtain a unified picture of the theory of abelian
and nonabelian class fields.

\subsection{Main Theorem of the Paper, An Introduction}

Given two integral algebraic schemes $X/Y$ over a fixed field $k$.
An algebraic $k-$scheme $X'$ is said to be a conjugate of $X$ over
$Y$ if there exists a $k-$isomorphism from $X$ onto $X'$ over $Y$,
and then $X$ is said to be Galois closed over $Y$ if there exists
one and only one conjugate of $X$ over $Y$. (See \S 1 for detail).

Assume that $\phi:X\rightarrow Y$ is a morphism of finite type.
Denote by $${Aut}_{k}\left( X/Y\right)$$ the group of
$k-$automorphisms of $X$ over $Y$. Let $\phi^{\sharp}(k(Y))$ be the
subfield of $k(X)$ induced from $\phi$ with Galois group
$$Gal(k\left( X\right) /\phi^{\sharp}(k(Y))).$$ Here $k(X)$ denotes the
field of rational functions on $X$. Now we are ready to relate the
main theorem of the present paper. (See \S 2 for detail).

\emph{\textbf{Theorem 2.1. (Main Theorem).} Let $X$ and $Y$ be two
integral $k-$varieties. Suppose that $X$ is Galois $k-$closed over
$Y$ by a surjective morphism $\phi$ of finite type.}

\emph{Then $k\left( X\right) $ is a Galois extension of
$\phi^{\sharp}(k(Y))$ and there is a group isomorphism
$${Aut}_{k}\left( X/Y\right) \cong Gal(k\left( X\right)
/\phi^{\sharp}(k(Y)))$$ if $k\left( X\right) /\phi^{\sharp}(k(Y))$ is separably generated.}

We will prove Theorem 2.1 in \S 3. From Theorem 2.1 it is seen that
Galois closed schemes behave like Galois extensions of fields, where
the groups of automorphism can be regarded as the Galois groups of
the extensions.

For example, take a number field $\mathbb{Q}\left( \xi\right)$. Let
$\xi^{\prime}$ be a $\mathbb{Q}-$conjugate of $\xi$. It is seen that
$Spec\left(\mathbb{Q}\left( \xi^{\prime}\right) \right) $ is a
conjugate of $Spec\left(\mathbb{Q}\left( \xi\right) \right) $.
Moreover,  let $\mathbb{Q}\left( \xi\right)/\mathbb{Q}$ be a Galois
extension. Then $Spec\left(\mathbb{Q}\left( \xi\right) \right) $ has
a unique conjugate, and hence is Galois closed, where the Galois
group is exactly the group of automorphisms.

We try to use the data of irreducible $k-$varieties $X/Y$ to
describe a given Galois extension $E/F$. We say that $X/Y$ are a
model for $E/F$ if the Galois group $Gal\left( E/F\right) $ is
isomorphic to the group of automorphisms of $X$ over $Y$.

The theorem above shows us some evidence that there can exist a nice
relationship between a Galois closed scheme and a Galois extension
of a field especially concerned with the nonabelian theory of class
fields.

\subsection{Outline of the Proof for the Main Theorem}

The whole of \S 3 will devote to prove the main theorem of the
present paper after we make definitions and fix notation in \S\S
1-2.

Here our approach to the proof will be established upon a full
analysis on affine open subsets of a given scheme with a preferable
favor of differential topology.

In \S 3.1 we will define conjugations of a given field. For a field
extension, the notion \textquotedblleft complete\textquotedblright\
is exactly the counterpart to the notion \textquotedblleft
normal\textquotedblright\ for an algebraic extension; the
\textquotedblleft conjugation\textquotedblright\ of a field is
exactly the counterpart to \textquotedblleft
conjugate\textquotedblright\ of a field for an algebraic extension.
Some results are proved by using the theory of specializations in a
scheme$^{[1]}$.

In \S 3.2 we will define conjugations of an open subset of a given
scheme in an evident manner. An open subset of a scheme is said to
have a complete set of conjugations if its conjugations all can be
affinely realized in the scheme. In deed, the definition here is the
geometric counterpart to the algebraic one in \S3.1.

Then we will establish a relationship between the conjugations of a
fixed field and the conjugations of an open subset of a given
scheme. So the discussions on fields and schemes are parallel.

In \S 3.3 we will prove that any finitely separably generated
extension is a Galois extension if and only it is complete, which is
equivalent to say that it contains its conjugations all. This is
very similar to an algebraic extension. Here, Weil's theory of
specializations$^{[9]}$ serves to prove the theorem.

To gain such results for schemes, we will be required to obtain
further properties for affine open sets in a given scheme. Thus in
\S 3.4 we will discuss affine structures on a scheme. Here there are
no new essential results and the discussion on affine structures is
just some interpretation of that in $[2]$. As usual, a scheme is a
ringed space covered by a family of affine schemes. By an affine
covering of a scheme we understand such a family of affine schemes.
Each affine covering of a scheme determines a unique affine
structure on the scheme. In general, a scheme can have many affine
structures on it; a given affine structure on a scheme can be
contained in many schemes, i.e., distinct structure sheaves on the
same underlying space, and such schemes are all isomorphic to each
other. However, for a Galois closed scheme, all affine structures on
it are contained in a unique scheme. This is one of the key points
to prove the main theorem of the paper.

Together with these preparations, in \S3.5 we will prove that each
affine open subset has a complete set of conjugations in integral
algebraic schemes which are Galois closed. It follows that the
fields of rational functions of such Galois closed schemes are
Galois extensions under surjective structure morphisms if the
extensions are separably generated. This gives a part of the proof
for the main theorem of the paper.

Finally in \S 3.6 we will complete the proof for the main theorem of
the paper. By isomorphisms of schemes, we will construct a
homomorphism $t$ between the group of automorphisms of the Galois
closed schemes and the Galois group of the fields of rational
functions on these schemes. It well-known that there exists a
bijection from the set of homomorphisms of algebras onto the set of
morphisms of their spectra. Then it is easily seen that the
homomorphism $t$ is injective. We will show any element of the
Galois group determines local isomorphisms on affine open subsets of
the scheme. All such local isomorphisms will patch an automorphism
of the whole of the scheme. This proves that the homomorphism $t$ is
surjective.

As the conclusion of this subsection, it should be noticed that
affine open subset sets in schemes work very well for the ramified
extensions of the fields of the rational functions on the schemes
while the points in the schemes may be good for unramified
extensions if one use the data of schemes to describe the theory
class fields.

\subsection*{Convention}

By a $k$\textbf{$-$variety} we will understand a scheme of finite
type over $Spec\left( k\right) $. We will follow throughout the
terminology of Grothendieck's EGA, except when otherwise specified.

\subsection*{Acknowledgment}

The author would like to express his sincere gratitude to Professor
Li Banghe for his invaluable advice and instructions on algebraic
geometry and topology.

\bigskip

\bigskip

\bigskip

\section{Definition for Galois Closed Schemes}

\bigskip

The notion of the conjugates of a $k-$variety is a generalization
from the classical affine varieties to algebraic schemes. Roughly
speaking, the conjugates of a $k-$variety behave as the conjugates
of a field, and a Galois closed $k-$variety behaves as a Galois
extension of a field.

\subsection{Galois Closed Schemes}

Fixed a field $k$ and a $k-$variety $S$. Let $X$ and $Y$ be two
$k-$varieties over $S$ by morphisms $f$ and $g$ respectively.

Then $X$ and $Y$ are \textbf{jointly of finite type} over $S$ if
there is an affine open covering $\{W_{\alpha }\} $ of $S$
satisfying the conditions:

$(i)$ Both $f^{-1}\left( W_{a}\right) $ and $g^{-1}\left( W_{\alpha
}\right) $ are finite unions of affine open sets $U_{\alpha i}$ and
$V_{\alpha j}$, respectively.

$(ii)$ Both $\mathcal{O}_{X}\left( U_{\alpha i}\right) $ and
$\mathcal{O}_{Y}\left( V_{\alpha j}\right) $ are algebras of finite
types over $\mathcal{O}_{S}\left( W_{a}\right)$.

It is immediate that two schemes are jointly of finite type over a
fixed scheme if and only if they are both of finite types over it
respectively.

Let $K/k$ be an extension and let $X$ and $Y$ be of finite type over
$S$. If there is a $K-$isomorphism $\sigma :X\rightarrow Y$ over
$S$, then $Y$ is a \textbf{$K-$conjugate} of $X$ over $S$, and
$\sigma $ is a \textbf{$K-$transformation} of $X$ onto $Y$ over $S$
. If $X=Y$, $\sigma $ is a \textbf{$K-$automorphism} of $X$ over
$S$. Put
\begin{equation*}
\begin{array}{l}
Conj_{K}\left( X/S\right)=\{K-\text{conjugates of
}X\text{ over }S\};\\

Aut_{K}\left( X/S\right)= \{K-\text{automorphisms of }X\text{ over
}S\}.
\end{array}
\end{equation*}

\begin{definition}
\emph{A $k-$variety $X$ over $S$ is said to be \textbf{Galois
$K-$closed} over $S$ if the identity
$(X,\mathcal{O}_{X})=(Y,\mathcal{O}_{Y})$ holds for any
$K-$conjugate $Y$ of $X$ over $S$.}
\end{definition}

\begin{remark}
\emph{There exists a nice relationship between Galois closed schemes
and Galois extensions of fields. Let $X$ be a Galois $k-$closed
scheme over a $k-$variety $S.$ Then $X/S$ can be intuitively
regarded as \textquotedblleft a Galois extension of the
field\textquotedblright\ with \textquotedblleft Galois
group\textquotedblright\ ${Aut}_{k}\left( X/S\right) .$}
\end{remark}

For example, $Spec\left( \mathbb{Q}\left( \sqrt{2}\right) \right) $
is Galois $\mathbb{Q}-$closed over $Spec\left(\mathbb{Q}\right)$
with a group isomorphism
$$
Aut_{\mathbb{Q}}\left( Spec\left( \mathbb{Q}\left( \sqrt{2}\right) \right)
\right) \cong Gal\left( \mathbb{Q}\left( \sqrt{2}\right) /\mathbb{Q}\right).
$$ It is immediate that $Spec\left( \mathbb{Q}\left(
^{3}\sqrt{2}\right) \right) $ is not Galois $\mathbb{Q}-$closed over
$Spec(\mathbb{Q})$ since $\mathbb{Q}\left( ^{3}\sqrt{2}\right) $ is
not a Galois extension.

\bigskip

\subsection{Galois Closures}

By a \textbf{Galois} \textbf{$K-$closure} of a $k-$variety $X$ over
$S$, denoted by $\overline{X}^{K}$, we understand a Galois
$K-$closed $k-$variety over $X$ which is a closed subscheme of any
other Galois $K-$closed $k-$variety over $X$.

It is immediate that $X$ is Galois $K-$closed over $S$ if and only
if there is $X=\overline{X}^{K}.$

There are some approaches to find a concrete Galois closure
$\overline{X}^{K}$ of a $k-$variety $X$. The finite group
actions$^{[3]}$ can afford us such a closure.

Here are some preliminary facts. Let $X/S$ be $k-$varieties such
that $Conj_{k}\left( X/S\right) $ is a finite set. Then there is a
$k-$variety $Y$ which is Galois $k-$closed over $S$ with $X$ a
closed subscheme of $Y$, and there is $\dim \overline{X}^{k}=\dim X$
if $\dim X<\infty$. Hence, any Artinian $k-$variety $X$ has a Galois
$k-$closure $\overline{X} ^{k}$ over $S$.

\bigskip

\bigskip

\section{Statement of The Main Theorem}

\bigskip

In the following we will relate the main result of the paper. It
shows us some evidence that there exists a nice relationship between
a Galois closed scheme and a Galois extension of a field especially
concerned with the nonabelian theory of class fields.

\subsection{Notation}

Let $D$ be an integral domain. Denote by $Fr(D)$ the field of
fractions on $D$. Given an extension $E$ of a field $F$ with Galois
group $Gal(E/F)$. Recall that $E/F$ is a Galois extension if $F$ is
the invariant subfield in $E$ for the Galois group $Gal(E/F)$. Here
$E/F$ is not necessarily an algebraic extension.

Let $X/Y$ be two irreducible $k-$varieties with morphism $\phi:X
\rightarrow Y$. Denote by $Aut_{k}(X/Y)$ the group of
$k-$automorphisms of $X$ over $Y$. Let $\xi$ be the generic point of
$X$. Then $\phi^{\sharp}(\mathcal{O}_{Y,\phi(\xi)})\subseteq
\mathcal{O}_{X,\xi}$ is a subring. Define
$$k(X)=Fr(\mathcal{O}_{X,\xi})$$
and
$$\phi^{\sharp}(k(Y))=Fr(\phi^{\sharp}(\mathcal{O}_{Y,\phi(\xi)})).$$

\bigskip

\subsection{Statement of the Main Theorem}

Here is the main theorem of the present paper.

\begin{theorem}
\textbf{(Main Theorem).} Let $X$ and $Y$ be two integral
$k-$varieties. Suppose that $X$ is Galois $k-$closed over $Y$ by a
surjective morphism $\phi$ of finite type.

Then $k\left( X\right) $ is a Galois extension of
$\phi^{\sharp}(k(Y))$ and there is a group isomorphism
$${Aut}_{k}\left( X/Y\right) \cong Gal(k\left( X\right)
/\phi^{\sharp}(k(Y)))$$ if $k\left( X\right) /\phi^{\sharp}(k(Y))$ is separably generated.
\end{theorem}

We will prove Theorem 2.1 in \S 3.

\begin{remark}
\emph{Theorem 2.1 affords us some certain information of the Galois
extension of a field in terms of data of the groups of the rational
automorphisms of schemes.}
\end{remark}

We attempt to use the data of schemes $X/Y$ to describe the field
extension $ E/F$, especially the (nonabelian) Galois group
$Gal(E/F)$. In particular, two irreducible $k-$varieties $X/Y$ are a
\textbf{$k-$model} of the field extension $E/F$ if there is a group
isomorphism
$$
Gal\left( E/F\right) \cong {Aut}_{k}\left( X/Y\right) $$
by a surjective morphism $f:X\rightarrow Y$.

We believe that in virtue of such data relating to schemes we can
obtain a unified picture of Galois extensions of fields from
algebraic ones to transcendental ones and from abelian ones to
nonabelian ones, where there exists the Galois correspondence and
class fields which are represented by affine open subsets of the
schemes.

\bigskip

\bigskip

\section{Proof of the Main Theorem}

\bigskip

In the following we will proceed in several subsections to prove the
main theorem of the paper.

\bigskip

\subsection{Definition for Conjugations of a Field}

Denote by $Fr\left( D\right) $ the fractional field of an integral
domain $D$. Let $K/k$ be a field extension.

\begin{definition}
\emph{$K$ is said to be \textbf{$k- $complete} (or \textbf{complete}
over $k$) if every irreducible polynomial $f(X)\in F[X]$ which has a
root in $K$ factors completely in $K\left[ X\right] $ into linear
factors for any intermediate field $k\subseteq F\subseteq K$.}
\end{definition}

Let $D^{\prime }/D$ be integral domains. $D^{\prime }$ is said to be
\textbf{$D-$complete} (or \textbf{complete} over $D$) if $Fr\left(
D^{\prime }\right) $ is $Fr\left( D\right) -$complete.

Given a finitely generated extension $E/k$. The elements
$$w_{1},w_{2},\cdots ,w_{n}\in E$$ are said to be a
\textbf{$(r,n)-$nice $k-$basis}
 of $E$ (or simply, a \textbf{nice $k-$basis}) if the following conditions are satisfied:

$E$ is generated by $w_{1},w_{2},\cdots ,w_{n}\in E$ over $k$;

 $w_{1},w_{2},\cdots ,w_{r}$ constitute a transcendental basis of $E$ over $k$;

$w_{r+1},w_{r+2},\cdots ,w_{n}$ are linearly independent over
$k(w_{1},w_{2},\cdots ,w_{r})$. Here $0\leq r\leq n$.

\begin{definition}
\emph{Let $E$ and $F$ be finitely generated extensions of a given
field $k$. $F$ is said to be a \textbf{$k-$conjugation} of $E$ (or a
\textbf{conjugation} of $E$ over $k$) if $F$ is contained in the
algebraic closure of $E$ and there is a $(r,n)-$nice $k-$basis
$w_{1},w_{2},\cdots ,w_{n} $ of $E$ such that $F$ is a conjugate of
$E$ over the $k(w_{1},w_{2},\cdots ,w_{r})$.}
\end{definition}

\bigskip

\subsection{Definition for Conjugations of an Open Set}

The discussion in this subsection (cf Definition 3.5 below) is a
counterpart to that in \S 3.1 (cf Definition 3.2 above). We will
extend the context of the conjugation from a field to an open set in
a scheme.

Let us recall some preliminary facts about specializations in a
scheme which are useful for us to study integral schemes.

Given a scheme $X$ and two points $x,y \in X$. Then $y$ is said to
be a \textbf{specialization} of $x$ in $X$, denoted by $x\rightarrow
y$ in $X$, if $y$ is contained in the (topological) closure of the
set $\{x\}$ (See $[1]$ for detail).

Let $x$ be a point in an affine scheme $Spec(A)$. Denote by $j_{x}$
the prime ideal of the ring $A$ corresponding to the point $x$.

\begin{lemma}
\emph{Let $X$ be an integral scheme. Take any $x,y\in X$ such that
$x\rightarrow y$ in $ X.$ Then there is a canonical ring
monomorphism
\begin{equation*}
i_{x,y}:\mathcal{O}_{X,y}\rightarrow \mathcal{O}_{X,x}.
\end{equation*}}
\end{lemma}

\begin{proof}
As $x\rightarrow y$ holds in $X$, there is an affine open subset $U$
of $X$ containing $x$ and $y$ from Lemma 1.8$^{[1]}$. Put
$U=Spec\left( A\right) $ and we have $x\rightarrow y$ in $U$; then
$j_{x}\subseteq j_{y}$ holds in the ring $A$. Let $S=A\smallsetminus
j_{y}$ and $T=A\smallsetminus j_{x}.$ As $S\subseteq T$, there is a
canonical homomorphism
\begin{equation*}
\rho _{A}^{T,S}:S^{-1}A\rightarrow T^{-1}A
\end{equation*}%
of the fractional rings. It is seen that $\rho _{A}^{T,S}$ is
injective. As $S^{-1}A\cong \mathcal{O}_{X,y}$ and $T^{-1}A\cong
\mathcal{O}_{X,x}$ hold, we obtain the canonical ring monomorphism
\begin{equation*}
i_{x,y}:\mathcal{O}_{X,y}\rightarrow \mathcal{O}_{X,x}
\end{equation*}%
factored by $\rho _{A}^{T,S}.$
\end{proof}

The \textbf{length} of the specialization $x\rightarrow y$ in $X$,
denoted by $l(x,y)$, is defined to be the supremum among the
integers $n$ such that there exist a chain of specializations
$$x_{0}=x \rightarrow x_{1} \rightarrow x_{2}
\rightarrow \cdot\cdot\cdot \rightarrow x_{n}=y$$ in $X$
(See $[1]$ for detail).

\begin{proposition}
\emph{Let $X$ be an irreducible scheme of finite dimension. Then
each morphism $\sigma:X\longrightarrow X$ has a fixed point in $X$.
In particular, the generic point of $X$ is an invariant point of any
surjective morphism $\delta$ of $X$ onto $X$.}
\end{proposition}

\begin{proof}
We have $l\left( X\right) =\dim X$ from Remark 2.3$^{\left[ 1\right]
}$, where $l(X)$ is the length of the space $X$ that is defined to
be the supremum of the lengths of specializations in $X$. Let $\xi $
be the generic point of $X$.

 For any $x\in X$, there is $\xi \rightarrow x$ in $X$. Then we have $\sigma \left( \xi
\right) \rightarrow \sigma \left( x\right) $ in $X$ by Proposition
1.3$^{[1]}$ which says that every morphism of schemes preserves the
specializations.  In particular, we choose $x=\xi $ and then obtain
a chain of the specializations
$$
\xi \rightarrow \sigma \left( \xi \right) \rightarrow \sigma ^{2}\left( \xi
\right) \rightarrow \cdots \rightarrow \sigma ^{n}\left( \xi \right)
\rightarrow \cdots
$$
in $X$.

We must have $\sigma ^{n}\left( \xi \right) =\sigma ^{n+1}\left( \xi
\right) $ for some $n\in \mathbb{N}$ since $l\left( X\right) <\infty .$ This
proves $\sigma \left( \sigma ^{n}\left( \xi \right) \right) =\sigma
^{n}\left( \xi \right) $.

For any $x \in X$ we have $\delta(\xi) \rightarrow \delta(x)$ since
any morphism preserves specilazations$^{[1]}$. As $\delta$ is
surjective, we have some $x_{0} \in X$ such that $\xi
=\delta(x_{0})$; then $\delta(\xi) \rightarrow \delta(x_{0})$, and
hence $\xi \rightarrow \delta(\xi) \rightarrow \delta(x_{0})=\xi$
holds. This proves $\delta(\xi)=\xi$.
\end{proof}

Consider an integral scheme $X$. Let $x\in X$ and let $\xi$ be the
generic point of $X$. From Lemma 3.3 we have the canonical
embeddings
\begin{equation*}
{Fr}(\mathcal{O}_{X}(U))\subseteq {Fr}(\mathcal{O}_{X,x})=\mathcal{O}_{X,\xi
}=k(\xi )=k\left( X\right)
\end{equation*}%
for every open set $U$ of $X$ containing $x$.

Now Let $X$ and $Y$ be integral $k-$varieties, and let $\varphi
:X\rightarrow Y$ be a morphism of finite type. Fixed a point $y \in
\varphi(X)$ and an affine open subset $V \subseteq$ of $Y$ with $V
\cap \varphi(X) \neq \emptyset$. For any affine open subset $U
\subseteq \varphi^{-1}(V)$ of $X$, the restriction $${\varphi\mid
_{U}} :(U,\mathcal{O}_{X}\mid _{U}) \longrightarrow
(V,\mathcal{O}_{Y}\mid _{V})$$ is a morphism of the open subschemes;
it follows that
$${\varphi}^{\sharp}(\mathcal{O}_{Y}(V)) \subseteq \mathcal{O}_{X}(U)$$ is a subalgebra.
This leads us to obtain the following
definitions.

Let $U_{1},U_{2}\subseteq \varphi ^{-1}(V)$ be open sets in $X$.
Assume that $Fr(\mathcal{O}_{X}(U_{1}))$ is a conjugation of
$Fr(\mathcal{O}_{X}(U_{2}))$ over
$Fr(\varphi^{\sharp}(\mathcal{O}_{Y}(V))).$ Then $%
U_{1}$ is said to be a \textbf{$V-$conjugation} of $U_{2}$, and
$Fr(\mathcal{O}_{X}(U_{1}))$ is said to be \textbf{affinely
realized} in $X$ by $U_{1}$.

Let $x,x^{\prime }\in \varphi ^{-1}\left( y\right) $ be two points
in $X$. Suppose that $Fr\left( \mathcal{O}_{X,x^{\prime }}\right) $
is a conjugation of $Fr\left( \mathcal{O}_{X,x}\right) $ over
$\varphi^{\sharp}(Fr\left( \mathcal{O}_{Y,y}\right))$. Then
$x^{\prime }\in \varphi ^{-1}\left( x\right) $ is said to be a
\textbf{$y-$conjugation} of $x,$ and the conjugation $Fr\left(
\mathcal{O}_{X,x^{\prime }}\right) $ is said to be \textbf{affinely
realized} in $X$ by $x^{\prime }$.

\begin{definition}
\emph{An open set $U\subseteq \varphi ^{-1}(V)$ in $X$ is said to
have a \textbf{complete set of $V-$conjugations} in $X$ if each
conjugation in $k(X)$ of $Fr(\mathcal{O}_{X}(U))$ over
$Fr(\varphi^{\sharp}(\mathcal{O}_{Y}(V)))$ can be affinely realized
by an open set in $X$. If $k\left( X\right) $ is replaced by
$\overline{k\left( X\right) },$ such a complete set is said to be
\textbf{absolutely complete}.}

\emph{A point $x\in \varphi ^{-1}\left( y\right) $ in $X$ is said to
have a \textbf{complete set of $y-$conjugations} in $X$ if each
conjugation in $k(X) $ of $Fr\left( \mathcal{O}_{X,x}\right) $ over
$\varphi^{\sharp}(Fr\left( \mathcal{O}_{Y,y}\right))$ can be
affinely realized by a point in $X$. If $k\left( X\right) $ is
replaced by $\overline{k\left( X\right) },$ such a complete set is
said to be \textbf{absolutely complete}.}
\end{definition}

\bigskip

\subsection{Conjugations and Galois Extensions}

We have the following result for conjugations of fields.
\begin{theorem}
Let $K/k$ be a finitely generated extension. The following
statements are equivalent.

$\left( i\right) $ $K$ is a complete field over $k$.

$\left( ii\right) $ Fixed any $x\in K$ and any subfield $k\subseteq
k_{x}\subseteq K$. Then each $k_{x}-$conjugation in $ \overline{K}$
of $k_{x}\left( x\right) $ is contained in $K$.

$\left( iii\right) $ Each $k-$conjugation in $\overline{K}$ of $K$ is contained in
$K$.
\end{theorem}

\begin{proof}
$\left( i\right) \implies \left( ii\right) .$ Let $x\in K$ and $k\subseteq
k_{x}\subseteq K.$ If $x$ is a varible over $k_{x},$ $k_{x}\left( x\right) $ is the
unique $k-$conjugation in $\overline{K}$ of $k_{x}\left( x\right) .$ If $x$ is
algebraic over $k_{x},$ a $k_{x}-$conjugation of $k_{x}\left( x\right)$ which is
exactly a $k_{x}-$conjugate of $k_{x}\left( x\right)$ is contained in $K$ by the assumption
that $K$ is $k-$%
complete; then all $k_{x}-$conjugates in $\overline{K}$ of
$k_{x}\left( x\right)$ is contained in $K$.

$\left( ii\right) \implies \left( iii\right) .$ Hypothesize that there is a $%
k-$conjugation $H\subseteq \overline{K}$ of $K$ is not contained in $K,$ that is,
$H\setminus K$ is a nonempty set. Take any $x\in H\setminus K,$ and put
$\sigma \left( x\right) \in K,$ where $\sigma :H\rightarrow K$ is an isomorphism
over $k.$

From $\left( ii\right) $ it is seen that $%
k\left( x\right) \in H\subseteq \overline{K}$ that is a $k-$conjugation in $%
\overline{K}$ of $k\left( \sigma \left( x\right) \right)$  is in contained in $K$. In
particular, $x$ belongs to $K$, where we will obtain a contradiction. This proves that
every $k-$conjugation in $\overline{K}$ of $K$ is in $K.$

$\left( ii\right) \implies \left( i\right) .$ Let $F$ be a field such that $%
k\subseteq F\subseteq K,$ and let $f\left( X\right) $ be an
irreducible polynomial over $F.$ Suppose that $x\in K$ satisfies the
equation $f\left(
x\right) =0.$ As an $F-$conjugation is an $F-$conjugate, from $%
\left( ii\right) $ it is seen that every conjugate $z\in \overline{F}$ of $%
x $ over $F$ is contained in $K$; then $K$ is complete over $k.$

$\left( iii\right) \implies \left( ii\right) .$ Take any $x\in K$ and any field $F$
such that $k\subseteq F\subseteq K$. If $x$ is a varible over $F,$ $F\left( x\right)
$ is the unique $F-$conjugation in $\overline{K}$ of $F\left( x\right) $ itself, and
hence $F\left( x\right) $ is contained in $K.$

Now suppose that $x$ is
algebraic over $F.$ In the following we will prove that each $F-$conjugates in $%
\overline{K}$ of $x$ is contained in $K$.

Let $z\in \overline{K}$ be an $F-$conjugate of $x,$ and let
\begin{equation*}
\sigma_{x} :F\left( x\right) \rightarrow F\left( z\right) ,x\longmapsto z
\end{equation*}
be the isomorphism over $F$.

If $F=K,$ we have $\sigma_{x} =id_{K};$ then $%
z=x\in K$. From now on, we suppose $F\not=K.$

Assume that $v_{1},v_{2},\cdots ,v_{m}$ are a $(s,m)-$nice $F\left(
x\right) -$basis of $K$. As $v_{1}$ is a varible over $F\left(
x\right) ,$ by
the $F-$isomorphism $\sigma_{x}$ we obtain an isomorphism $\sigma _{1}$ of $%
F\left( x,v_{1}\right) $ onto $F\left( z,v_{1}\right) $ defined in an evident manner
that
\begin{equation*}
\sigma_{1}:\frac{f(v_{1})}{g(v_{1})}\mapsto \frac{\sigma _{x}(f)(v_{1})}{\sigma
_{x}(g)(v_{1})}
\end{equation*}%
for any polynomials
\begin{equation*}
f[X],g[X]\in F\left( x\right) [X]
\end{equation*}%
with $g[X]\neq 0$.

It is easily seen that $g(v_{1})=0$ if and only if ${\sigma
_{x}(g)(v_{1})}=0$. Hence, $\sigma_{1}$ is well-defined.

Similarly, for the varibles $v_{1},v_{2},\cdots ,v_{s}$ over
$F\left( x\right) ,$ there is a field isomorphism
$$
\sigma _{s}: F\left( x,v_{1},v_{2},\cdots ,v_{s}\right)\longrightarrow
F\left( z,v_{1},v_{2},\cdots ,v_{s}\right)
$$
defined by%
$$
x\longmapsto z\text{ and }v_{i}\longmapsto v_{i}
$$
for every $1\leq i\leq s.$ We have the restrictions $$\sigma_{i+1}|_{F\left(
x,v_{1},v_{2},\cdots ,v_{i}\right)}=\sigma_{i}$$ for $1\leq i\leq s-1$.

Consider $v_{s+1}.$ We have an isomorphism $\sigma _{s+1}$ of
$F\left( x,v_{1},v_{2},\cdots ,v_{s+1}\right) $ onto $F\left(
z,v_{1},v_{2},\cdots ,v_{s+1}\right) $ defined by
\begin{equation*}
x\longmapsto z\text{ and }v_{i}\longmapsto v_{i}
\end{equation*}
for every $1\leq i\leq s+1.$

Prove that $\sigma_{s+1}$ is well-defined. That is, we prove that
$f\left( v_{s+1}\right)  =0$ holds if and only if $\sigma_{s}\left(
f\right)  \left( v_{s+1}\right)  =0$ holds for any polynomial
$f\left(  X_{s+1}\right)  \in F\left(
x,v_{1},v_{2},\cdots,v_{s}\right)  \left[  X_{s+1}\right]  .$

It reduces to prove the following claim.

\qquad\textbf{Claim.} \emph{Given any $f\left(  X,X_{1},X_{2}
,\cdots,X_{s+1}\right)$ contained in the polynomial ring $F\left[
X,X_{1},X_{2},\cdots,X_{s+1}\right] .$ Then
 $f\left( x,v_{1},v_{2},\cdots,v_{s+1}\right)  =0$  holds if and only if $f\left(
z,v_{1},v_{2},\cdots,v_{s+1}\right)  =0$ holds.}

We use Weil's algebraic theory of specializations$^{\left[  9\right]  }$ to prove the
above claim. By Proposition 1 (Page 3 of $\left[  9\right]  $) it is seen that $F\left(
x\right)  $ and $F\left(  v_{1},v_{2},\cdots ,v_{s+1}\right)  $ are independent over
$F$ since $x$ is assumed to be algebraic over $F.$ From Weil's definition for
specializations it is clear that $\left(  z\right)  $ is a (generic) specialization of $\left(
x\right) $ over $F$. Then it follows that $\left(  z\right)  $ is a specialization of
$\left(  x\right)  $ over the field $F\left(  v_{1},v_{2},\cdots ,v_{s+1}\right)  $ in
virtue of Theorem 4 (Page 29 of $\left[  9\right]  $). This proves \textquotedblleft
only if\textquotedblright\ in the claim.

As $x$ and $z$ are conjugates over $F,$ it is seen that \textquotedblleft
if\textquotedblright\ in the claim is true if we substitute $z$ for $x.$

In such a manner we have a field isomorphism
\begin{equation*} \sigma _{m}: F\left(
x,v_{1},v_{2},\cdots ,v_{m}\right)\longrightarrow F\left( z,v_{1},v_{2},\cdots
,v_{m}\right)
\end{equation*}
defined by%
\begin{equation*}
x\longmapsto z\text{ and }v_{i}\longmapsto v_{i}
\end{equation*}
for every $1\leq i\leq m.$

Then $F\left( z,v_{1},v_{2},\cdots ,v_{m}\right)$ is a conjugation of $K=F\left(
x,v_{1},v_{2},\cdots ,v_{m}\right)$ over $F$. From the assumption $%
\left( iii\right) ,$ we have $F\left( z,v_{1},v_{2},\cdots ,v_{m}\right)\subseteq K.$
This proves $z\in K.$
\end{proof}

\bigskip

\subsection{Preliminaries: $k-$Affine Structures}

For the convenience of context, let us recall some preliminary
results on affine structures on an algebraic $k-$scheme in order to
obtain some further properties for Galois closed schemes. Here the
discussion about affine structures is just a word-by-word
interpretation of that in $[2]$ by substituting $k-$algebras for
rings, and there are no new essential results in this section.

By definition, a scheme $(X,\mathcal{O}_{X})$ is a locally ringed
space that can be covered by a family
$\{(U_{\alpha},\phi_{\alpha})\}_{\alpha \in \Delta}$ of affine
schemes. That is, for each $\alpha \in \Delta$ there is an
isomorphism $\phi_{\alpha}:U_{\alpha}\cong Spec{A_{\alpha}}$ with
$\{U_{\alpha}\}_{\alpha \in \Delta}$ an open covering of $X$ and
$A_{\alpha}$ a commutative ring of identity. We say that
$\{(U_{\alpha},\phi_{\alpha})\}_{\alpha \in \Delta}$ is an
\textbf{affine covering} of $(X,\mathcal{O}_{X})$. It is easily seen
that a scheme can have many affine coverings.

An affine covering $\{U_{\alpha},\phi_{\alpha}\}_{\alpha \in
\Delta}$ of $(X,\mathcal{O}_{X})$ is said to be \textbf{reduced} if
$U_{\alpha}\neq U_{\beta}$ holds for any $\alpha\neq \beta$ in
$\Delta$. We will denote by $(X,\mathcal{O}_{X};\mathcal{A}_{X})$ a
scheme $(X,\mathcal{O}_{X})$ together with a given reduced affine
covering $\mathcal{A}_{X}$.

Let $(X,\mathcal{O}_{X})$ and $(Y,\mathcal{O}_{Y})$ be two schemes
with reduced affine coverings $\mathcal{A}_{X}$ and
$\mathcal{A}_{Y}$ respectively. Then we say
$(X,\mathcal{O}_{X};\mathcal{A}_{X})=(Y,\mathcal{O}_{Y};\mathcal{A}_{Y})$
if and only if the following conditions are satisfied:

    $(i)$ As schemes we have $(X,\mathcal{O}_{X})=(Y,\mathcal{O}_{Y})$.

    $(ii)$ There is a reduced affine covering $\mathcal{A}$ of $(X,\mathcal{O}_{X})$ such
    that $\mathcal{A}_{X}$ and $\mathcal{A}_{Y}$ are both subsets of $\mathcal{A}$.

\bigskip

In the following we will give the further discussion on affine
coverings. From the discussion below it is easily seen that an
affine covering of a scheme determines a unique affine structure on
it.

Let $\mathfrak{Comm}{/k}$ be the category of finitely generated
algebras (with identities) over a given field $k$.

\begin{definition}
\emph{A \textbf{pseudogroup of }$k-$\textbf{affine transformations}, denoted by
$\Gamma$, is a set of isomorphisms of $k-$algebras satisfying the conditions (i)-(v):}

\emph{$\left(  i\right)  $ Each $\sigma\in\Gamma$ is a
$k-$isomorphism of $k-$algebras from $dom\left(  \sigma\right)  $
onto $rang\left( \sigma\right)$ contained in $\mathfrak{Comm}/k$.}

\emph{$\left(  ii\right)  $ If $\sigma\in\Gamma,$ the inverse
$\sigma^{-1}$ is contained in $\Gamma.$}

\emph{$\left(  iii\right)  $ The identity map $id_{A}$ on $A$ is
contained in $\Gamma$ for any $k-$algebras $A\in \mathfrak{Comm}/k$
if there is some $\delta\in\Gamma$ with $dom\left(  \delta\right)
=A.$}

\emph{$\left(  iv\right)  $ If $\sigma\in\Gamma,$ the
$k-$isomorphism induced by $\sigma$ defined on the localization
$dom\left(  \sigma\right)  _{f}$ at any $f\in dom\left(
\sigma\right)  $ is contained in $\Gamma.$}

\emph{$\left(  v\right)  $ Let $\sigma,\delta\in\Gamma.$ The
$k-$isomorphism factorized by $dom\left(  \tau\right)  $ from
$dom\left( \sigma\right)  _{f}$ onto $rang\left(  \delta\right)
_{g}$ is contained in $\Gamma$ if $\tau\in\Gamma$ holds and there
are $k-$isomorphisms $ dom\left(  \tau\right)  \cong dom\left(
\sigma\right)  _{f}$ and
$
dom\left(  \tau\right)  \cong rang\left(  \delta\right)  _{g}%
$
for some $f\in dom\left(  \sigma\right)  $ and $g\in rang\left(
\delta\right)  .$}
\end{definition}

Let $X$ be a topological space, and $\Gamma$ a pseudogroup of $k-$affine
transformations.

\begin{definition}
\emph{A $k-$\textbf{affine }$\Gamma-$\textbf{atlas}
$\mathcal{A}\left( X,\Gamma\right)  $ on $X$ is a collection of pairs $\left(
U_{j},\varphi _{j}\right) ,$ called $k-$\textbf{affine charts}, satisfying the
conditions (i)-(iii):}

\emph{$\left(  i\right)  $ For every $\left(
U_{j},\varphi_{j}\right)  \in \mathcal{A}\left(  X,\Gamma\right)  ,$
$U_{j}$ is an open subset of $X$ and $\varphi _{j}$ is an
homeomorphism of $U_{j}$ onto $Spec\left(  A_{j}\right)  ,$ where
$A_{j}$ is a $k-$algebra contained in $\Gamma$.}

\emph{$\left(  ii\right)  $ $\bigcup U_{j}$ is an open covering of
$X.$}

\emph{$\left(  iii\right)  $ Given any $\left(
U_{i},\varphi_{i}\right) ,\left( U_{j},\varphi_{j}\right)  \in
\mathcal{A}\left(  X,\Gamma\right)  $ with $U_{i}\cap U_{j}\not
=\varnothing$. There exists $\left(  W_{ij},\varphi _{ij}\right) \in
\mathcal{A}\left( X,\Gamma\right)  $ such that $W_{ij}\subseteq
U_{i}\cap U_{j}$, and the $k-$isomorphism from the localization
$\left( A_{j}\right)  _{f_{j}}$ onto the localization $\left(
A_{i}\right)  _{f_{i}}$  which is induced by the restriction
$$
\varphi_{j}\circ\varphi_{i}^{-1}\mid_{W_{ij}}:\varphi_{i}(W_{ij}%
)\rightarrow\varphi_{j}(W_{ij})
$$
is contained in $\Gamma$. Here $A_{i}$ and $A_{j}$ are $k-$algebras
contained in $\Gamma$ such that $ \varphi_{i}\left(  U_{i}\right)
=SpecA_{i}$ and $\varphi_{j}\left( U_{j}\right)  =SpecA_{j}$ hold
and there are homeomorphisms
$
\varphi_{i}\left(  W_{ij}\right)  \cong Spec\left(  A_{i}\right)  _{f_{i}%
}$ and $\varphi_{j}\left(  W_{ij}\right)  \cong Spec\left(
A_{j}\right)
_{f_{j}}%
$
for some $f_{i}\in A_{i}$ \text{ and }$f_{j}\in A_{j}.$}
\end{definition}

A $k-$affine $\Gamma-$atlas $\mathcal{A}\left(
X,\Gamma\right) $ on $X$ is said to be \textbf{complete} (or \textbf{maximal}) if it
can not be contained properly in any other $k-$affine $\Gamma-$atlas of $X.$

\begin{definition}
\emph{Two $k-$affine $\Gamma-$atlases $\mathcal{A}$ and
$\mathcal{A}^{\prime}$ on $X$\ are said to be $\Gamma-$\textbf{compatible} if
the following condition is satisfied:}

\emph{Given any $\left(  U,\varphi\right)  \in \mathcal{A}$ and
$\left( U^{\prime
},\varphi^{\prime}\right)  \in \mathcal{A}^{\prime}$ with $U\cap U^{\prime}%
\not =\varnothing.$ There exists a $k-$affine chart $\left(
W,\varphi ^{\prime\prime}\right)  \in \mathcal{A}\bigcap
\mathcal{A}^{\prime}$ such that $W\subseteq U\cap U^{\prime}$, and
the $k-$isomorphism from the localization $ A _{f}$ onto the
localization $\left(  A^{\prime}\right)  _{f^{\prime}}$ induced by
the restriction $\varphi^{\prime}\circ\varphi^{-1}\mid_{W}$ is
contained in $\Gamma$, where $A$ and $A^{\prime}$ are $k-$algebras
contained in $\Gamma$ such that $ \varphi\left( U\right)  =SpecA$
and $\varphi^{\prime}\left(  U^{\prime
}\right)  =SpecA^{\prime}%
$ hold, and there are homeomorphisms $ \varphi\left(  W\right)
\cong SpecA _{f}$ and $ \varphi^{\prime}\left(  W\right)  \cong
Spec\left( A^{\prime}\right) _{f^{\prime}}$ for some $f\in A$ and
$f^{\prime}\in A^{\prime}.$}
\end{definition}

Let $X$ be a topological space. By a $k-$\textbf{affine }$\Gamma-$%
\textbf{structure} on $X$ we understand a complete $k-$affine $\Gamma-$atlas
$\mathcal{A}\left( \Gamma\right)  $ on $X$ for some pseudogroup $\Gamma$ of
$k-$affine transformations.

Fixed a pseudogroup $\Gamma$ of $k-$affine
transformations. By Zorn's Lemma it is seen that for any given $k-$affine $\Gamma-$atlas
$\mathcal{A}$ on $X$ there is a unique complete $k-$affine $\Gamma-$atlas
$\mathcal{A}_{m}$ on $X$ such that

$\left(  i\right)  $ $\mathcal{A}\subseteq \mathcal{A}_{m};$

$\left(  ii\right)  $ $\mathcal{A}$ and $\mathcal{A}_{m}$ are
$\Gamma-$compatible.

Here $\mathcal{A}$ is said to be a \textbf{base} for $\mathcal{A}_{m}$ and
$\mathcal{A}_{m}$ is the complete $k-$affine $\Gamma-$atlas determined by
$\mathcal{A}.$

\begin{definition}
\emph{Given a $k-$affine $\Gamma-$structure $\mathcal{A}\left(
\Gamma\right)  $ on the space $X$. Assume that there exists a
locally ringed space $\left(  X,\mathcal{F}\right)  $ such that for
each $\left(  U_{\alpha },\varphi_{\alpha}\right)  \in
\mathcal{A}\left(  \Gamma\right)  $ there is $
\varphi_{\alpha\ast}\mathcal{F}\mid_{U_{\alpha}}\left(
SpecA_{\alpha}\right) =A_{\alpha} $, where $A_{\alpha}$ is a
$k-$algebra contained in $\Gamma$ with $\varphi _{\alpha}\left(
U_{\alpha}\right)  =SpecA_{\alpha}.$}

\emph{Then $\mathcal{A}\left(  \Gamma\right)  $ is said to be an \textbf{admissible }%
$k-$\textbf{affine structure} on $X,$ and $\left(  X,\mathcal{F}\right)  $ is
an \textbf{extension} of the $k-$affine $\Gamma-$structure $\mathcal{A}\left(
\Gamma\right)  .$}
\end{definition}

\begin{remark}
\emph{All extensions of an admissible affine structure on a
topological space are schemes which are isomorphic with each
other$^{[2]}$.}
\end{remark}

\bigskip

Now take an algebraic $k-$scheme $\left(  X,\mathcal{O}_{X}\right) $
with a fixed affine covering $\mathcal{C}_{X}$.

Denote by $\Gamma_{0}$ (respectively, $\Gamma^{max}$)  the
union of the set of some (respectively, all) identities of $k-$algebras
$
id_{A_{\alpha}}:A_{\alpha}\rightarrow A_{\alpha}%
$ and the set of some (respectively, all) isomorphisms of
$k-$algebras $ \sigma_{\alpha\beta}:\left(  A_{\alpha}\right)
_{f_{\alpha}}\rightarrow \left(  A_{\beta}\right)  _{f_{\beta}}, $
satisfying the conditions $(i)-(ii)$:

$(i)$ Each $A_{\alpha},A_{\beta},A_{\gamma}\in Comm$ are
$k-$algebras such that there are affine open subsets
$
U_{\alpha},U_{\beta},$ and $U_{\gamma}\subseteq U_{\alpha}\cap U_{\beta}%
$ of $X$ with identities $ \varphi_{\alpha}\left(  U_{\alpha}\right)
=SpecA_{\alpha},$ $\varphi_{\beta}\left(  U_{\beta}\right)
=SpecA_{\beta},$ and $ \varphi_{\gamma}\left(  U_{\gamma}\right)
=SpecA_{\gamma}.
$

$(ii)$ Each $\sigma_{\alpha\beta}:\left(  A_{\alpha}\right)  _{f_{\alpha}%
}\rightarrow\left(  A_{\beta}\right)  _{f_{\beta}}$ is induced from
the homeomorphism $
\varphi_{\alpha}\circ\varphi_{\beta}^{-1}\mid_{U_{\gamma}}:\varphi_{\beta
}({U_{\gamma}})\rightarrow\varphi_{\alpha}({U_{\gamma}}) $ of spaces
such that $ \varphi_{\alpha}\left(  U_{\gamma}\right)  \cong
Spec\left(  A_{\alpha }\right)  _{f_{\alpha}}$ and
$\varphi_{\beta}\left(  U_{\gamma
}\right)  \cong Spec\left(  A_{\beta}\right)  _{f_{\beta}}%
$
hold for some $f_{\alpha}\in A_{\alpha}$ and $f_{\beta}\in A_{\beta}.$

Then the pseudogroup generated by $\Gamma_{0}$ in $\mathfrak{Comm}/k$,
which is defined to be the smallest pseudogroup
containing $\Gamma_{0}$ in $\mathfrak{Comm}/k$, is called a \textbf{pseudogroup of $k-$affine
transformations in }$\left(  X,\mathcal{O}_{X}\right)  $ (relative to the given affine covering).

In particular, the pseudogroup generated by $\Gamma^{max}$ in $\mathfrak{Comm}/k$ is
called the \textbf{maximal pseudogroup
of $k-$affine transformations in }$\left(  X,\mathcal{O}_{X}\right)  $ (relative
to the given affine covering).

\begin{definition}
\emph{Let $\Gamma$ be a pseudogroup of $k-$affine transformations  in an algebraic $k-$scheme
$\left(  X,\mathcal{O}_{X}\right)  $. Define
$$
\mathcal{A}^{\ast}\left(  \Gamma\right)  =\{\left(
U_{\alpha},\varphi_{\alpha}\right)  :\varphi_{\alpha}\left(  U_{\alpha
}\right)  =SpecA_{\alpha}\emph{ and }A_{\alpha}\in\Gamma\}
$$
where each $U_{\alpha}$ is an affine open subset in the scheme $X$.}

\emph{If $\mathcal{A}^{\ast}\left(
\Gamma\right)  $ is a $k-$affine $\Gamma-$atlas on the underlying space $X$, then $\Gamma$ is said to be a
\textbf{canonical pseudogroup of $k-$affine transformations} in
$\left(  X,\mathcal{O}_{X}\right)  $, and $\mathcal{A}^{\ast}\left(  \Gamma
\right)  $ is called a \textbf{$k-$affine atlas in
}$(X,\mathcal{O}_{X})$.}
\end{definition}

\begin{definition}
\emph{Let $\Gamma$ be a canonical pseudogroup of $k-$affine
transformations in an algebraic $k-$scheme $\left(
X,\mathcal{O}_{X}\right) $, and let $\mathcal{A}$ be a $k-$affine
$\Gamma-$atlas on the underlying space $X$.}

\emph{$(i)$ $\mathcal{A}$
is said to be a
\textbf{canonical $k-$affine structure} in $(X,\mathcal{O}_{X})$ if
$\mathcal{A}^{\ast}\left(  \Gamma\right) $ is a base for $\mathcal{A}$.}

\emph{$(ii)$ $\mathcal{A}$ is said to be a \textbf{relative
canonical $k-$affine structure} in $(X,\mathcal{O}_{X})$ if $\mathcal{A}$ is maximal among
all the $k-$affine $\Gamma-$atlases in $(X,\mathcal{O}_{X})$ which contain
$\mathcal{A}^{\ast}\left(  \Gamma\right)  $ and are $\Gamma-$compatible with
$\mathcal{A}^{\ast}\left(  \Gamma\right)  $.}
\end{definition}

It is seen that the relative
canonical $k-$affine structure is well-defined.
In deed, the $k-$affine $\Gamma-$atlases in
$(X,\mathcal{O}_{X})$ which are $\Gamma-$compatible with
$\mathcal{A}^{\ast}\left(  \Gamma\right)  $
are all are $\Gamma-$compatible with each other. The converse is true if
they are assumed to contain
$\mathcal{A}^{\ast}\left(  \Gamma\right)  $.

Fixed such a canonical pseudogroup $\Gamma$ in $(X,\mathcal{O}_{X})$.
By Zorn's Lemma it is easily seen that such a canonical affine $\Gamma-$structure in $\left(
X,\mathcal{O}_{X}\right)  $ is unique. However,
there can be many canonical affine structure in $\left(
X,\mathcal{O}_{X}\right)  $ when $\Gamma$ varies.

$\left(  X,\mathcal{O}_{X}\right)$ is said to have \textbf{a
unique} (respectively, \textbf{relative}) \textbf{canonical
$k-$affine structure} if there exists only one (respectively,
relative) canonical $k-$affine structure in it.

Evidently, any affine open subset $U$ in $X$  is contained in a
canonical $k-$affine $\Gamma -$structure if and only if it is
contained in a relative canonical $k-$affine $\Gamma-$structure. It
follows that an algebraic $k-$scheme has a unique canonical
$k-$affine structure if and only if it has a unique relative one.

\begin{remark}
\emph{Given an algebraic $k-$scheme $\left(
X,\mathcal{O}_{X}\right)$.}

$(i)$ \emph{Let $\Gamma$ be the maximal pseudogroup of $k-$affine
transformations in $\left(  X,\mathcal{O}_{X}\right) $. Then
$\mathcal{A}^{\ast}\left(  \Gamma\right)  $ is a relative canonical
$k-$affine $\Gamma-$structure, called an \textbf{intrinsic
$k-$affine structure} of $\left( X,\mathcal{O}_{X}\right)$. Denote
$\mathcal{A}^{\ast}\left( \Gamma\right)  $ by $\mathcal{A}_{X}$.}

$(ii)$ \emph{A given algebraic $k-$scheme can have many intrinsic
$k-$affine structures. Conversely, a given intrinsic $k-$affine
structures can be of many algebraic $k-$schemes.}

$(iii)$ \emph{An intrinsic $k-$affine structure of a given algebraic
$k-$scheme affords us the definition how the $k-$affine charts (ie,
the affine schemes) are patched into the scheme.}
\end{remark}

\begin{remark}
\emph{To be precise, an algebraic scheme $\left(
X,\mathcal{O}_{X}\right)$ should be defined by three types of data:}

    \qquad  \emph{$X$, the underlying space;}

    \qquad  \emph{$\mathcal{O}_{X}$, the sheaf on $X$;}

    \qquad  \emph{$\mathcal{A}$, an fixed intrinsic affine
    structure.}

\emph{Conversely, all these data $X$, $\mathcal{O}_{X}$, and
$\mathcal{A}$ completely determines a unique algebraic scheme,
denoted by $\left( X,\mathcal{O}_{X}; \mathcal{A}\right)$.}
\end{remark}

\begin{remark}
\emph{Let $\left( X,\mathcal{O}_{X}\right)$ be an algebraic
$k-$scheme defined in the usual manner. Given an affine covering $\mathcal{A}_{X}$
 of $\left( X,\mathcal{O}_{X}\right)$.
Then $\mathcal{A}_{X}$ is an atlas on the underlying space $X$, which
 determines a unique
intrinsic $k-$affine structure $\mathcal{A}$ by (i) of Remark
3.14 above. That is, $\mathcal{A}_{X}$ is a base for $\mathcal{A}$. In such a case,
we will identify $\left( X,\mathcal{O}_{X};
\mathcal{A}_{X}\right)$ with $\left( X,\mathcal{O}_{X};
\mathcal{A}\right)$.}
\end{remark}

\begin{definition}
\emph{An \textbf{associate scheme} of a given algebraic $k-$scheme
$\left(  X,\mathcal{O}_{X}\right)  $ is an extension on the
underlying space $X$ of a canonical affine structure or a relative
one in $\left(  X,\mathcal{O}_{X}\right)  $.}
\end{definition}

\begin{remark}
$(i)$ \emph{Any algebraic  $k-$scheme has an associate scheme. In
particular, it itself is an associate scheme of its own.}

$(ii)$ \emph{The associate schemes of a given algebraic  $k-$scheme
are all isomorphic with each other.}
\end{remark}

\bigskip

\subsection{Conjugations and Galois Closed Varieties}

The discussion in this subsection (cf Theorem 3.22 below) is a
counterpart to that in \S 3.3 (cf Theorem 3.6 above) as well.

\begin{proposition}
\emph{Let $X$ and $Y$ be integral $k-$varieties. Suppose that $X$ is Galois
$k-$closed over $Y$ by a morphism $\varphi$. Take any affine open subset $V$ of
$Y$ with
$\varphi(X)\bigcap V\neq\varnothing$.}

\emph{Then every affine open set
$U\subseteq\varphi^{-1}(V)$ of $X$ has a complete set of
$V-$conjugations in $X$, and such a complete set is absolutely
complete; moreover, $\mathcal{O}_{X}$ is the unique structure
sheaf on the underlying space $X$ with such a property.}
\end{proposition}

\begin{proof}
Let $U_{0}$ be an open subset of $X$. It is easily seen that $U_{0}$ is
 an affine open subset in $X$ if and only if there is an intrinsic $k-$affine
structure of $X$ containing $U_{0}$ as a $k-$affine chart. It
follows that there is some reduced affine covering of $X$ containing
$U_{0}$. Obviously, every reduced affine covering of $X$ determines
a unique intrinsic $k-$affine structure of $X$.

Now take any affine open set $V$ of $Y$ with $\varphi(X)\bigcap
V\neq\varnothing$. Let $U\subseteq\varphi^{-1}(V)$ be an affine open
subset of $X$.

$\left(  i\right)$ Prove that there is a reduced affine covering $\mathcal{A}_{\infty}$
of $X$ such that each conjugation of $U$ over $V$ is affinely
realized by some affine open set contained in
$\mathcal{A}_{\infty}$. We will proceed in two steps.

\textsl{Step 1.} Fixed any intrinsic $k-$affine structure of $X$ an
intrinsic $k-$affine structure $\mathcal{A}$ of $X$ containing $U$
as a $k-$affine chart. Then there exists a reduced affine covering
$\mathcal{A}_{X}$ of $X$ that determines $\mathcal{A}$.

Assume that $H$ is a conjugation of ${Fr}(\mathcal{O}_{X}(U))$ over
${Fr\left(\varphi^{\sharp}( \mathcal{O}_{Y}(V))\right)}$
such that $H$ can not be
affinely realized in $X$ by any affine open subset contained in
$\mathcal{A}_{X}$.

As $X$ is an integral scheme, there is an affine open subset $W$
contained in $U$ such that $$Fr(\mathcal{O}_{X}(W))\cong
Fr(\mathcal{O}_{X}(U)) \cong H.$$

Let $\sigma$ be an isomorphism of ${Fr}(\mathcal{O}_{X}(W))$ onto
$H$. Put $$B=\mathcal{O}_{X}(W) \text{ and } A=\sigma(B).$$ It is
seen that ${Fr}(A)=H$ holds.

We have isomorphisms
$$\varphi:\left(  W,\mathcal{O}_{W}\right)  \cong\left(  SpecB,\widetilde
{B}\right)  $$
and
$$\delta:\left(  SpecB,\widetilde{B}\right)  \cong\left(
SpecA,\widetilde{A}\right)  $$
where $\delta$ is induced from the isomorphism $\sigma:B\rightarrow A$ of $k-$algebras.

It follows that there is a reduced affine covering $\mathcal{B}_{X}$
of $X$ such that
$$\mathcal{B}_{X}=\mathcal{A}_{X}\bigcup \{(W_{H},\psi_{H})\}$$
with $(W_{H},\psi_{H})=(W,\delta \circ \phi)$.

Hence, for each conjugation $H$ of ${Fr}(\mathcal{O}_{X}(U))$ over
${Fr}\left(\varphi^{\sharp}( \mathcal{O}_{Y}(V))\right)$ there is a reduced affine
covering $\mathcal{B}_{H}$ such that $H$ can be affinely realized by
an affine open set $U_{H}$ contained in $\mathcal{B}_{H}$. In
particular, each $U_{H}$ is contained in $\psi^{-1}(V)$.

Moreover, given any two conjugations $H_{1}\neq H_{2}$ of
${Fr}(\mathcal{O}_{X}(U))$ over ${Fr}\left(\varphi^{\sharp}( \mathcal{O}_{Y}(V))\right)$.
There are affine open sets
$U_{H_{1}}\neq U_{H_{2}}$ contained in some reduced affine covering
$\mathcal{B}_{H_{1,2}}$ of $X$ such that $H_{1}$ and ${H_{2}}$ can
be affinely realized by $U_{H_{1}}$ and $U_{H_{2}}$ respectively.

\textsl{Step 2.} Fixed a reduced affine covering
$\mathcal{A}_{X}=\{(U_{\alpha},\phi_{\alpha})\}$ of $X$.

Let
$$\mathcal{B}^{*}=\bigcup_{H}\{(U_{H},\phi_{H})\}$$ where $H$ runs
through all conjugations of ${Fr}(\mathcal{O}_{X}(U))$ over the field
${Fr}\left(\varphi^{\sharp}( \mathcal{O}_{Y}(V))\right)$, and $(U_{H},\phi_{H})$ is
contained in a reduced affine covering of $X$ such that $H$ is
affinely realized in $X$ by $U_{H}$.

Put
$$\mathcal{A}_{X}^{*}=\{(U_{\alpha},\phi_{\alpha}):U_{\alpha}=U_{H}\text{ holds for some }
 (U_{H},\phi_{H}) \in  \mathcal{B}^{*}\}.$$

Then we obtain a reduced affine covering
$$\mathcal{A}_{\infty}=(\mathcal{A}_{X} \setminus \mathcal{A}_{X}^{*}) \bigcup \mathcal{B}^{*}$$
such that every conjugation of $U$ over $V$ is affinely realized in
$X$ by some affine open set contained in $\mathcal{A}_{\infty}$.

This proves that the affine open set $U\subseteq\varphi^{-1}(V)$ of
$X$ has a complete set of $V-$conjugations in $X$.

By replacing $k\left(  X\right)  $ by $\overline{k\left(  X\right)
},$ it is seen that such a complete set is absolutely complete. In
deed, it automatically holds that such a set is absolutely complete
for integral schemes.

$\left(  ii\right)  $ Now prove the uniqueness of the structure sheaf $\mathcal{O}_{X}$.
Let $\mathcal{B}_{\infty}$ be the $k-$affine structure on the
underlying space $X$ determined by $\mathcal{A}_{\infty}$.

Then $\mathcal{B}_{\infty}$ is admissible on the space $X$ since the
sheaf $\mathcal{O}_{X}$ is an extension of $\mathcal{B}.$

Take any extension $\mathcal{F}$ of $\mathcal{B}_{\infty}$ on the
space $X.$ From Remark 3.11 it is seen that the scheme $\left(
X,\mathcal{F}\right)  $ is isomorphic to the scheme $\left(
X,\mathcal{O}_{X}\right)$, and it follows that
$\left(  X,\mathcal{F}\right)  $ is a $k-$conjugate of $\left(  X,\mathcal{O}%
_{X}\right)  $ over $Y$.

Hence, we have $\left( X,\mathcal{F}\right) =\left(
X,\mathcal{O}_{X}\right)$ since $X$ is Galois $k-$closed over $Y$.
This proves $\mathcal{F}=\mathcal{O}_{X}$.
\end{proof}

\begin{theorem}
Let $X$ and $Y$ be integral $k-$varieties and let $X$ be Galois
$k-$closed over $Y$ by a surjective morphism $\varphi $ of finite
type.

Then each conjugation in $\overline{k\left(
X\right) }$ of $k\left( X\right) $ over $\varphi^{\sharp}(k(Y))$ is contained in $k\left( X\right) $.
\end{theorem}

\begin{proof}
Hypothesize there is a conjugation $H$ in $\overline{%
k\left( X\right) }$ of $k\left( X\right) $ over
$\varphi^{\sharp}(k(Y))$ is not contained in $k\left( X\right)$.
Suppose $u_{0}\in H\setminus k\left( X\right)$.

Let $\sigma :H\rightarrow k\left( X\right) $ be an isomorphism over
$\varphi^{\sharp}(k(Y))$. Put $$ v_{0}=\sigma \left( u_{0}\right)
.$$

Prove that there are affine open subsets $V$ of $Y$ and $U\subseteq
\varphi ^{-1}\left( V\right) $ of $X$ such that $v_{0}\in
\mathcal{O}_{X}\left( U\right) $.

In fact, let $\xi $ and $\eta $ be the generic points of $X$ and
$Y$, respectively. We have $$\varphi (\xi )=\eta,
\mathcal{O}_{X,\xi }={Fr}(\mathcal{O}_{X,\xi }),\text{ and }\mathcal{O}%
_{Y,\eta }={Fr}(\mathcal{O}_{Y,\eta }). $$

Then
$$
v_{0}\in \mathcal{O}_{X,\xi }$$ and $$\mathcal{O}%
_{X,\xi }=\lim_{\overrightarrow{W}}\mathcal{O}_{X}\left( W\right)$$
hold, where $W$ runs through all open sets in $X$. It follows that
there is some open $W_{1}$ of $X$ such that $v_{0}$ belongs to
$\mathcal{O}_{X}( W_{1}) $.

It is clear that there are affine open subsets $V_{\alpha}$ (with
$\alpha \in \Gamma$) of $Y$ such that $$\bigcup_{\alpha \in
\Gamma}V_{\alpha}\supseteq \varphi \left( W_{1}\right) .$$

Let $V_{\alpha_{0}}$ be the affine open subset such that
$$V_{\alpha_{0}}\bigcap \varphi(W_{1}) \neq \emptyset.$$

Then
$$\varphi^{-1} (V_{\alpha_{0}})\bigcap W_{1} \neq \emptyset ,$$
which is an open subset of $X$.

Now take an affine open subset $W_{\alpha_{0}}$ of $X$ such that
$$W_{\alpha_{0}} \subseteq \varphi^{-1} (V_{\alpha_{0}})\bigcap W_{1}.$$

It is seen that $v_{0}$ is contained in
$\mathcal{O}_{X}(W_{\alpha_{0}})$ from the given injective
homomorphism $$\mathcal{O}_{X}(W_{1})\longrightarrow
\mathcal{O}_{X}(W_{\alpha_{0}}).$$

Hence, we obtain affine open subsets $$V=V_{\alpha_{0}}$$ of $Y$ and
$$U=W_{\alpha_{0}}\subseteq \varphi ^{-1}\left( V\right) $$ of $X$
such that $$v_{0}\in \mathcal{O}_{X}\left( U\right) .$$

By Proposition 3.19 it is seen that $U$ has an absolutely complete set of $%
V-$conjugations in $X$ and that $\mathcal{O}_{X}$ is the unique
structure sheaf with that property, and it follows that $H$ can be
affinely realized in $X$ by an affine open set $U^{\prime }$. Then we
have $$u_{0}\in H=Fr\left( \mathcal{O} _{X}\left( U^{\prime }\right)
\right)\cong  Fr\left( \mathcal{O} _{X}\left( U\right) \right) \ni
v_{0}.$$

Hence, $u_{0}$ is contained in $k\left( X\right)$, which is in
contradiction to the above hypothesis that $u_{0} \notin k(X)$. This
completes the proof.
\end{proof}

\begin{corollary}
\emph{Let $X$ and $Y$ be integral $k-$varieties and let $X$ be
Galois $k-$closed over $Y$ by a surjective morphism $\varphi $. Then
$k\left( X\right) $ is a complete extension over
$\varphi^{\sharp}(k\left( Y\right)) $.}
\end{corollary}

\begin{proof}
It is immediate from Theorems 3.6 and 3.20.
\end{proof}

\bigskip

\subsection{Proof of The Main Theorem}

Now we are ready to prove the main theorem of the paper, Theorem 2.1
in \S 2.

\begin{proof}
It is immediate that $k\left( X\right) /\phi^{\sharp}(k(Y))$ is a
Galois extension by Theorem 3.6 and Corollary 3.21.

It is easily seen that $k(X)$ is the set of all the elements of the
forms
$$(U,f)$$ with $U$ an open set of $X$ and $f$ an element of $\mathcal{O}_{X}(U)$.
That is, $k(X)$ is the field of rational functions on $X$. Here we
identify $\mathcal{O}_{X}(U)$ with its image in $k(X)$ since the
homomorphism $$\mathcal{O}_{X}(U) \longrightarrow k(X)$$ of rings is
injective for every open subset $U$ of $X$. Let $\xi$ be the generic
point of $x$.

In the following we prove that there is a group isomorphism
$${Aut}_{k}\left( X/Y\right) \cong Gal(k\left( X\right)
/\phi^{\sharp}(k(Y))).$$
We will proceed in several steps.

\bigskip

\emph{Step 1.} Take any
$$
\sigma =\left( \sigma ,\sigma ^{\sharp}\right) \in Aut_{k}\left( X/Y\right) .
$$ That is, $$\sigma : X \longrightarrow X$$ is a homeomorphism, and $$\sigma ^
{\sharp}:\mathcal{O}_{X} \rightarrow \sigma _{\ast
}\mathcal{O}_{X}$$ is an isomorphism of sheaves of rings on $X$. It
follows that
$$
\sigma ^{\sharp}:k\left( X\right)=\mathcal{O}_{X,\xi } \rightarrow \sigma _{\ast
}\mathcal{O}_{X,\xi }=k\left( X\right)
$$
is an automorphism of $k(X)$. Let $ \sigma ^{\sharp-1}$ be the
inverse of $\sigma ^{\sharp}$.

Fixed any open subset $U$ of $X$. We have the restriction
$$\sigma=(\sigma ,\sigma ^{\sharp}): (U,\mathcal{O}_{X}|_{U}) \longrightarrow
(\sigma(U),\mathcal{O}_{X}|_{\sigma(U)})$$ of open subschemes. In
particular, $$\sigma^{\sharp}:\mathcal{O}_{X}|_{\sigma(U)}
\rightarrow \sigma_{\ast}\mathcal{O}_{X}|_{\sigma(U)}$$ is an
isomorphism of sheaves. For every $f \in \mathcal{O}_{X}|_{U}(U)$,
there is $$f \in
\sigma_{\ast}\mathcal{O}_{X}|_{\sigma(U)}(\sigma(U)),$$ and hence
$$\sigma^{\sharp -1}(f) \in \mathcal{O}_{X}(\sigma(U)).$$

Define a mapping
$$
t:Aut_{k}\left( X/Y\right) \longrightarrow Gal\left( k\left(
X\right) /\phi^{\sharp}(k\left( Y\right) )\right)$$ of sets by
$$\sigma \longmapsto t(\sigma)=\left\langle \sigma ,\sigma ^{\sharp
-1}\right\rangle
$$
such that
$$
\left\langle \sigma ,\sigma ^{\sharp-1}\right\rangle :\left( U,f\right)
 \longmapsto \left( \sigma \left( U\right) ,\sigma
^{\sharp-1}\left( f\right) \right) $$ is the mapping of $k(X)$ into
$k(X)$ induced by $\sigma \in Aut_{k}\left( X/Y\right)$.

\emph{Step 2.} Prove that $t$ is well-defined. In deed, given any
$$
\sigma =\left( \sigma ,\sigma ^{\sharp}\right) \in Aut_{k}\left( X/Y\right).
$$ For any $(U,f),(V,g) \in k(X)$, we have
$$(U,f)+(V,g)=(U\cap V, f+g)$$ and
$$(U,f)\cdot(V,g)=(U\cap V, f\cdot g).$$

Then we have
\begin{equation*}
\begin{array}{l}
\left\langle \sigma ,\sigma ^{\sharp-1}\right\rangle((U,f)+(V,g))\\

=\left\langle \sigma ,\sigma ^{\sharp-1}\right\rangle((U\cap V,
f+g))\\

=(\sigma(U\cap V),
\sigma^{\sharp -1}(f+g))\\

=(\sigma(U\cap V), \sigma^{\sharp -1}(f))+(\sigma(U\cap V),
\sigma^{\sharp -1}(g))\\

=(\sigma(U), \sigma^{\sharp -1}(f))+(\sigma(V), \sigma^{\sharp
-1}(g))\\

=\left\langle \sigma ,\sigma ^{\sharp-1}\right\rangle((U,f))+
\left\langle \sigma ,\sigma ^{\sharp-1}\right\rangle((V,g)).
\end{array}
\end{equation*}
and
\begin{equation*}
\begin{array}{l}
\left\langle \sigma ,\sigma ^{\sharp-1}\right\rangle((U,f)\cdot(V,g))\\

=\left\langle \sigma ,\sigma ^{\sharp-1}\right\rangle((U\cap V,
f\cdot g))\\

=(\sigma(U\cap V),
\sigma^{\sharp -1}(f\cdot g))\\

=(\sigma(U\cap V), \sigma^{\sharp -1}(f))\cdot(\sigma(U\cap V),
\sigma^{\sharp -1}(g))\\

=(\sigma(U), \sigma^{\sharp -1}(f))\cdot(\sigma(V), \sigma^{\sharp
-1}(g))\\

=\left\langle \sigma ,\sigma ^{\sharp-1}\right\rangle((U,f))\cdot
\left\langle \sigma ,\sigma ^{\sharp-1}\right\rangle((V,g)).
\end{array}
\end{equation*}

It follows that $ \left\langle \sigma ,\sigma
^{\sharp-1}\right\rangle$ is an automorphism of $k\left( X\right) .$

Consider the given morphism
$\phi=(\phi,\phi^{\sharp}):(X,\mathcal{O}_{X})\rightarrow
(Y,\mathcal{O}_{Y})$ of schemes. By Proposition 3.4 it is seen that
$\phi(\xi)$ is the generic point of $Y$ and that $\xi$ is an
invariant point of each automorphism $\sigma \in Aut_{k}\left(
X/Y\right)$. Then $\sigma^{\sharp}: \mathcal{O}_{X,\xi} \rightarrow
\mathcal{O}_{X,\xi}$ is an isomorphism of rings.

On the other hand, it is seen that $\sigma^{\sharp}$ is an
isomorphism over $\phi^{\sharp}(k(Y))$ since $\mathcal{O}_{X,\xi}$
is an algebra over $\phi^{\sharp}(\mathcal{O}_{Y,\phi(\xi)})$ and
$\sigma^{\sharp}$ is an automorphism of $\mathcal{O}_{X,\xi}$ over
$\phi^{\sharp}(\mathcal{O}_{Y,\phi(\xi)})$.

Hence,
$$
\left\langle \sigma ,\sigma ^{\sharp-1}\right\rangle \in Gal\left( k\left(
X\right) /\phi^{\sharp}(k\left( Y\right) \right)) .
$$

Now take any
$$
\sigma =\left( \sigma ,\sigma ^{\sharp}\right) ,\delta =\left( \delta ,\delta
^{\sharp}\right) \in Aut_{k}\left( X/Y\right) .
$$
We have
$$
\left\langle \delta ,\delta ^{\sharp-1}\right\rangle \circ \left\langle \sigma
,\sigma ^{\sharp-1}\right\rangle =\left\langle \delta \circ \sigma ,\delta
^{\sharp-1}\circ \sigma ^{\sharp-1}\right\rangle
$$
since
$$
\delta ^{\sharp -1}\circ \sigma^{\sharp -1}=(\delta \circ \sigma)^{\sharp -1}
$$ holds.

Hence,
$$
t:Aut_{k}\left( X/Y\right) \rightarrow Gal\left( k\left( X\right) /\phi^{\sharp}\left(k\left(
Y\right)\right) \right)
$$
is a homomorphism of groups.

\emph{Step 3.} Prove that ${t}$ is surjective. In fact, given any
element $\rho$ of the group $Gal\left( k\left( X\right)
/\phi^{\sharp}\left(k\left( Y\right)\right) \right) $. We have
$$
\rho :\left( U_{f},f\right) \in k\left( X\right) \longmapsto \left( U_{\rho
\left( f\right) },\rho \left( f\right) \right) \in k\left( X\right) ,
$$ where $U_{f}$ and $U_{\rho (f)}$ are open subsets in $X$ such that $f \in \mathcal{O}_{X}(U_{f})$
and $\rho (f) \in \mathcal{O}_{X}(U_{\rho (f)})$.

Fixed any $k-$affine structure $\mathcal{A}$ on $X$. Let $\left(
U,\varphi \right)\in \mathcal{A}$ be a $k-$affine chart with
$\varphi \left( U\right) =SpecA_{U}$. We have
$$
A_{U}\cong \mathcal{O}_{X}(U) =\{\left( U_{f},f\right) \in
k\left( X\right) :U_{f}\supseteq U\}.
$$

Put
$$
B_{U}=\{\left( U_{\rho \left( f\right) },\rho \left( f\right) \right) \in
k\left( X\right) :U_{\rho(f)}\supseteq U\}.
$$

We have $B_{U} \subseteq \mathcal{O}_{X}(U) $. As $\rho$ is
surjective, each element $(W_{g},g)$ in $\mathcal{O}_{X}(U)$ is the
image of some element $(V_{h},h)$ in $k(X)$ under $\rho$; then
$B_{U} \supseteq \mathcal{O}_{X}(U) $, and hence $B_{U} =
\mathcal{O}_{X}(U) $. This proves $$A_{U} \cong B_{U}$$ and
$$\rho (\mathcal{O}_{X}(U))=\mathcal{O}_{X}(\rho(U)).$$

It is seen that there is a unique isomorphism
$$\lambda_{U}=\left(\lambda_{U}, \lambda_{U}^{\sharp} \right): (U, \mathcal{O}_{X}|_{U})
\rightarrow (U, \mathcal{O}_{X}|_{\rho(U)})$$ of the affine open subscheme in $X$ such that
$$\rho |_{\mathcal{O}_{X}(U)}=\lambda_{U}^{\sharp -1}: \mathcal{O}_{X}(U) \rightarrow \mathcal{O}_{X}(\rho(U)).$$

Now we show that there is an automorphism $\lambda$ of scheme $X$
such that $$\lambda|_{U}=\lambda_{U}$$ holds for each affine open
subscheme $U$ of $X$.

In fact, take any affine open subsets $U$ and $V$ of $X$. As
morphisms of schemes, it is seen that $$\lambda_{U}|_{U\cap
V}=\lambda_{V}|_{U\cap V}$$ holds since we have
$$\rho |_{\mathcal{O}_{X}(U\cap V)}=\lambda_{U}|_{U\cap
V}^{\sharp -1}: \mathcal{O}_{X}(U\cap V) \rightarrow
\mathcal{O}_{X}(\rho(U\cap V))$$ and
$$\rho |_{\mathcal{O}_{X}(U\cap V)}=\lambda_{V}|_{U\cap
V}^{\sharp -1}: \mathcal{O}_{X}(U\cap V) \rightarrow
\mathcal{O}_{X}(\rho(U\cap V))$$ by the above construction for each $\lambda_{U}$.

Let $U$ and $V$ be any affine open subsets of $X$. It is seen that
for any points $x,y\in X$ there is
$$\lambda_{U}(x)=\lambda_{V}(y)$$ if and only if $$x,y\in U\cap V $$
holds. In deed, take an affine open subset $W$ of $X$ with
$z=\lambda_{U}(x) \in W$. If $x$ can not be contained in $V$, we
will have affine open subsets $x\in U_{0}\subseteq U$, $y\in
V_{0}\subseteq V$, and $z\in W_{0}\subseteq W$ which are isomorphic
to each other as schemes such that
$$\lambda_{U}(U_{0})=W_{0}=\lambda_{V}(V_{0})$$ and that $$V_{0}\not\ni x\in
U_{0}, y\in V_{0},\text{ and }z\in W_{0};$$ then there will be
$f_{0}\in\mathcal{O}_{X}(U_{0})$, $g_{0}\in\mathcal{O}_{X}(V_{0})$,
and $h_{0}\in\mathcal{O}_{X}(W_{0})$ such that $$\rho
\left((U_{0},f_{0})\right)=(W_{0},h_{0})=\rho
\left((V_{0},g_{0})\right)$$ but $$(U_{0},f_{0})
\not=(V_{0},g_{0}),$$ which is in contradiction to the assumption
that $\rho$ is an isomorphism.

Then we have a homeomorphism $\lambda$ of $X$ onto $X$ as a
topological space defined in an evident manner that
 $$\lambda: x\in X \mapsto
\lambda_{U}(x)\in X$$ if $x$ is contained in an affine open subset
$U$ of $X$. The mapping $\lambda$ is well-defined since all affine
open subsets of $X$ constitute a base for the topology on $X$.
Hence, we obtain an isomorphism $\lambda\in Aut_{k}\left( X\right)
.$

We show that $\lambda\in Aut_{k}\left( X/Y\right)$ holds with
$t\left(\lambda\right)=\rho$. In deed, as $\rho$ is an isomorphism
of $k(X)$ over $\phi^{\sharp}\left(k(Y)\right)$, it is seen that the
isomorphism $\lambda_{U}$ is over $Y$ by $\phi$ for any affine open
subset $U$ of $X$; then $\lambda$ is an automorphism of $X$ over $Y$
by $\phi$. It is immediate that $t\left(\lambda\right)=\rho$ holds.

This proves that there exists $\lambda\in Aut_{k}\left( X/Y\right) $
such that $t(\lambda)=\rho$ for each $\rho \in Gal\left( k\left(
X\right) /\phi^{\sharp}\left(k\left( Y\right)\right) \right) $. So,
$ {t}$ is a surjection.

\emph{Step 4.} Prove that ${t}$ is injective. Assume $\sigma ,\sigma
^{\prime }\in {Aut}_{k}\left( X/Y\right) $ such that $t\left( \sigma
\right) =t\left( \sigma ^{\prime }\right) .$ There is
$$
\left( \sigma \left( U\right) ,\sigma ^{\sharp-1}\left( f\right) \right) =\left(
\sigma ^{\prime }\left( U\right) ,\sigma ^{\prime \sharp-1}\left( f\right)
\right)
$$
for any $\left( U,f\right) \in k\left( X\right) .$ In particular, we
have
$$
\left( \sigma \left( U_{0}\right) ,\sigma ^{\sharp-1}\left( f\right) \right)
=\left( \sigma ^{\prime }\left( U_{0}\right) ,\sigma ^{\prime \sharp-1}\left(
f\right) \right)
$$
for each affine open subset $U_{0}$ of $X$ such that
$$
\sigma \left( U_{0}\right) =\sigma ^{\prime }\left( U_{0}\right) \subseteq
\sigma \left( U\right) \cap \sigma ^{\prime }\left( U\right)
$$
with $f\in \mathcal{O}_{X}(U_{0})$; then $\sigma |_{U_{0}}=\sigma
^{\prime }|_{U_{0}}$ holds as isomorphisms of schemes. Hence, we
have $ \sigma = \sigma ^{\prime }$. This proves that $t$ is an
injection.

At last we obtain an isomorphism
$$
t:{Aut}_{k}\left( X/Y\right) \cong Gal\left( k\left( X\right)
/\phi^{\sharp}\left(k\left(
Y\right)\right) \right)
$$ of groups.
This completes the proof.
\end{proof}

\end{document}